\documentclass[amsthm]{elsart}
\makeatletter
\def\ps@copyright{\let\@mkboth\@gobbletwo%
 \let\@oddhead\@empty
 \let\@evenhead\@empty
 \def\@oddfoot{\scriptsize\fontsize{8}{13}\selectfont\slshape\hskip-0em
   To be published in Journal of Applied Probability Vol.\ 52 No.\ 2 \copyright2015 by The Applied Probability Trust}%
 \let\@evenfoot\@oddfoot}
\makeatother
\usepackage{graphicx}
\usepackage{amssymb}
\usepackage[final]{microtype}
\usepackage{amsmath}
\usepackage{enumerate,xcolor}
\theoremstyle{plain}
\newtheorem{Sat}{Definition}[section]
\newtheorem{The}[Sat]{Theorem}
\newtheorem{Bem}[Sat]{Remark}
\newcommand{\Var}{\operatorname{Var}}
\newcommand{\Erw}{\operatorname{E}}
\newcommand{\Prob}{\operatorname{P}}
\begin{document}
\begin{frontmatter}

\title{On the Acceleration of the Multi-Level Monte Carlo Method}

\author{Kristian Debrabant\thanksref{label1}} and
\ead{debrabant@imada.sdu.dk}
\author{Andreas R\"{o}{\ss}ler\thanksref{label2}}
\ead{roessler@math.uni-luebeck.de}
\address[label1]{University of Southern Denmark, Department of Mathematics and Computer Science (IMADA),
Campusvej 55, 5230 Odense M, Denmark}
\address[label2]{Universit\"{a}t~zu~L\"ubeck, Institut~f\"ur~Mathematik, Ratzeburger~Allee~160,
23562~L\"ubeck, Germany}
\begin{abstract}
The multi-level Monte Carlo method proposed by M.~Giles~(2008)
approximates the expectation of some functionals applied to a
stochastic process with optimal order of convergence for the
mean-square error. In this paper, a modified multi-level Monte Carlo
estimator is proposed with significantly reduced computational
costs. As the main result, it is proved that the modified estimator
reduces the computational costs asymptotically by a factor
$(p/\alpha)^2$ if weak approximation methods of orders $\alpha$ and
$p$ are applied in case of computational costs growing with same
order as variances decay.
\end{abstract}

\end{frontmatter}
\allowdisplaybreaks
\section{Introduction} \label{Introduction}
%
%
The multi-level Monte Carlo method proposed in \cite{Gil08a}
approximates the expectation of some functional applied to some
stochastic processes like e.~g.\ solutions of stochastic
differential equations (SDEs) at a lower computational complexity
than classical Monte Carlo simulation, see also
\cite{Gil08b,Hei01,Ke05}. Multi-level Monte Carlo approximation is
applied in many fields like mathematical finance \cite{Av09,Gil09a},
for SDEs driven by a L\'{e}vy process \cite{Der11}, by fractional
Brownian motion \cite{KNP11} or for stochastic PDEs \cite{SchGit11}.
The main idea of this article is to
reduce the computational costs additionally by applying the multi-level
Monte Carlo method as a variance reduction technique for some higher
order weak approximation method. As a result, the computational
effort can be significantly reduced while the optimal order of
convergence for the root mean-square error is preserved.

The outline of this paper is as follows. We give a brief
introduction to the main ideas and results of the multi-level Monte
Carlo method in Section~\ref{Section2:MLMC-Simulation-Original}.
Based on these results, in
Section~\ref{Sec3:Improved-MLMC-Estimator} we present as the main
result a modified multi-level Monte Carlo algorithm that allows to
reduce the computational costs significantly. Depending on the
relationship between the orders of variance reduction and of the
growth of the costs, there exists a reduction of the computational
costs by a factor depending on the weak order of the underlying
numerical method. As an example, the modified multi-level Monte
Carlo algorithm is applied to the problem of weak approximation for
stochastic differential equations driven by Brownian motion in
Section~\ref{Sec4:Numerical-Examples-SDEs}.
%
%
\section{Multi-level Monte Carlo simulation}
\label{Section2:MLMC-Simulation-Original}
%
%
Let $(\Omega, \mathcal{F}, \Prob)$ be a probability space with some
filtration $(\mathcal{F}_t)_{t \geq 0}$ and
let $X=(X_t)_{t \in I}$
denote an adapted stochastic process on the interval $I=[t_0,T]$
that belongs to a space $\mathbb{X}$ that may be infinite dimensional.
In the following, we are interested in the approximation of
$\Erw_{\Prob}(f(X))$ for some functional $f \in \mathbb{F}$ where
$\mathbb{F}$ denotes a suitable class of functionals that are of
interest.
Further, let an equidistant discretization $I_h = \{t_0, t_1,
\ldots, t_N\}$ with $0 \leq t_0 < t_1 < \ldots < t_N =T$ of the time
interval $I$ with step size $h$ be given. Then, we consider a
probability space $(\tilde{\Omega}, \tilde{\mathcal{F}},
\tilde{\Prob})$ with some filtration $(\tilde{\mathcal{F}}_t)_{t \in
I_h}$ and we denote by $Y=(Y_{t})_{t \in I_h}$ a discrete time
approximation of $X$ on the grid $I_h$, adapted to
$(\tilde{\mathcal{F}}_t)_{t \in I_h}$.
Thus, we consider the
approximation $Y \in \mathbb{X}_h$ of $X \in \mathbb{X}$ on a finite dimensional space
$\mathbb{X}_h$.
Here, the probability spaces
$(\Omega, \mathcal{F}, \Prob)$ and $(\tilde{\Omega},
\tilde{\mathcal{F}}, \tilde{\Prob})$ may be but do not have to be
equal and we assume that $Y$ approximates $X$ in the weak sense with
some order $p >0$, i.e.\
\begin{equation}
    | \Erw_{\tilde{\Prob}}(f(Y)) - \Erw_{\Prob}(f(X)) | = {O}(h^p)
\end{equation}
for all $f \in \mathbb{F}$.

In order to approximate the expectation of $f(X)$ we apply the
multi-level Monte Carlo estimator introduced in \cite{Gil08a}. For
some fixed $M \in \mathbb{N}$ with $M \geq 2$ and some $L \in
\mathbb{N}$ we define the step sizes $h_l = \frac{T}{M^{l}}$ and let
$Y^l=(Y_t)_{t \in I_{h_l}}$ denote the discrete time approximation
process on the grid $I_{h_l}$ based on step size $h_l$ for $l=0,1,
\ldots, L$.
Here, we consider the
approximations $Y^l \in \mathbb{X}_{h_l}$ for $l=0,1, \ldots, L$ of $X \in \mathbb{X}$ on a
sequence $\mathbb{X}_{h_0} \subset \mathbb{X}_{h_1} \subset \ldots \subset \mathbb{X}_{h_L}$
of finite dimensional subspaces.
Then, the multi-level Monte Carlo estimator is defined by
\begin{equation} \label{MLMC-estimator-Giles}
    \hat{Y}_{ML} = \sum_{l=0}^L \hat{Y}^l
\end{equation}
for some $L \in \mathbb{N}$ using the estimators $\hat{Y}^0 =
\frac{1}{N_0} \sum_{i=1}^{N_0} f({Y^0}^{(i)})$ and
\begin{equation}
    \hat{Y}^l = \frac{1}{N_l} \sum_{i=1}^{N_l} \left( f({Y^l}^{(i)})
    - f({Y^{l-1}}^{(i)}) \right)
\end{equation}
for $l=1, \ldots, L$.
Then, we get
\begin{equation}
    \begin{split}
    \Erw_{\tilde{\Prob}}(\hat{Y}_{ML})
    &= \Erw_{\tilde{\Prob}}(f(Y^0)) + \sum_{l=1}^L \Erw_{\tilde{\Prob}}(f(Y^l) - f(Y^{l-1})) \, .
    \end{split}
\end{equation}
Here, we have to point out, that both approximations ${Y^l}^{(i)}$
and ${Y^{l-1}}^{(i)}$ are simulated simultaneously based on the same
realisation of the underlying driving random process whereas
$({Y^l}^{(i)},{Y^{l-1}}^{(i)})$ and $({Y^l}^{(j)},{Y^{l-1}}^{(j)})$
are independent realisations for $i \neq j$.

Now, there are two sources of errors for the approximation. On the
one hand, we have a systematical error
that depends on the dimension of $\mathbb{X}_{h_l}$
due to the discrete time
approximation $Y^l \in \mathbb{X}_{h_l}$ based on step size $h_l$
which is given by the bias of the method. On the other hand, there is a statistical error
from the estimator for the expectation of $f(Y^l)$ by the Monte Carlo
simulation. Therefore, we consider the root mean-square error
\begin{equation} \label{Sec2:root-mean-square-error}
    \begin{split}
    e(\hat{Y}_{ML}) &= \left( \Erw_{\tilde{\Prob}}(| \hat{Y}_{ML}
    - \Erw_{\Prob}(f(X)) |^2) \right)^{1/2}
    \end{split}
\end{equation}
of the multi-level Monte Carlo method in the following. In order to
rate the performance of an approximation method, we will analyse the
root mean-square error of the method compared to the computational
costs. Therefore, we denote by $C(Y)$ the computational costs of the
approximation method $Y$. In order to determine $C(Y)$, one may use
a cost model where e.g.\ each operation or evaluation of some
function is charged with the price of one unit, i.e.\ one counts the
number of needed mathematical operations or function evaluations.
Further, each random number that has to be generated to compute $Y$
may also be charged with the price of one unit.

It is well known that the optimal order of convergence
for the classical Monte Carlo estimator $\hat{Y}_{MC} =
\frac{1}{N} \sum_{i=1}^N f({Y}^{(i)})$ is given by \[e(\hat{Y}_{MC})
= O \left( (1/C(\hat{Y}_{MC}) )^{\frac{p}{2p+1}} \right)\] where $p$
is the weak order of convergence of the approximations $Y$, see
Duffie and Glynn~\cite{DuGl95}. Thus, higher order weak
approximation methods result in a higher order of convergence with
respect to the root mean-square error. Clearly, the best root
mean-square order of convergence that can be achieved is at most
$1/2$. However, the order bound $1/2$ can not be reached by any weak
order $p$ approximation method in the case of the classical Monte
Carlo simulation. Therefore, in order to attain the optimal order of
convergence for the root mean-square error we apply the multi-level
Monte Carlo estimator (\ref{MLMC-estimator-Giles}).
The following theorem due to Giles~\cite{Gil08a} is presented in a
slightly generalized version suitable for our considerations.
\begin{The} \label{Main-Theorem-Giles}
    For some $L \in \mathbb{N}$, let $Y^l$
    denote the approximation process on the grid $I_{h_l}$ with respect
    to step size $h_l=\frac{T}{M^l}$ for each $l=0,1, \ldots, L$, respectively.
    Suppose that there exist some constants $\alpha > 0$ and
    $c_{1,\alpha}, c_{2,0}, c_{2}, c_{2,L} >0$
    and $\beta,\beta_L>0$ such that for the bias
    \begin{enumerate}
        \item[1)] $|\Erw_{\Prob}(f(X))-\Erw_{\tilde{\Prob}}(f(Y^{L}))| \leq c_{1,\alpha} \, h_L^{\alpha}$
    \end{enumerate}
    and for the variances
    \begin{enumerate}
        \item[2)] $\Var_{\tilde{\Prob}}(f(Y^0)) \leq c_{2,0} \, h_0^{\beta}$,
        \item[3)] $\Var_{\tilde{\Prob}}(f(Y^l)-f(Y^{l-1})) \leq c_{2} \, h_l^{\beta}$ \quad for $l=1, \ldots, L-1$,
        \item[4)] $\Var_{\tilde{\Prob}}(f(Y^L)-f(Y^{L-1})) \leq c_{2,L} \, h_L^{\beta_{L}}$.
    \end{enumerate}
    Further, assume that there exist constants $c_{3,0}, c_3, c_{3,L} >0$
    and $\gamma, \gamma_{L} \geq 1$ such that for the computational costs
    \begin{enumerate}
        \item[5)] $C(Y^0) \leq c_{3,0} \, T \, h_0^{-\gamma}$,
        \item[6)] $C(Y^l,Y^{l-1}) \leq c_{3} \, T \, h_l^{-\gamma}$ \quad for $l=1, \ldots, L-1$,
        \item[7)] $C(Y^L,Y^{L-1}) \leq c_{3,L} \, T \, h_L^{-\gamma_{L}}$.
    \end{enumerate}
    Then, for some arbitrarily prescribed error bound $\varepsilon>0$ there exist values $L$
    and $N_l$ for $l=0,1, \ldots, L$, such that the root mean-square error of the
    multi-level Monte Carlo estimator $\hat{Y}_{ML}$ has the bound
    \begin{equation}
        e(\hat{Y}_{ML}) < \varepsilon
    \end{equation}
    with computational costs bounded by
    \begin{equation}
        C(\hat{Y}_{ML}) \leq \begin{cases}
                              c_4 \, \varepsilon^{-2} & \text{ if } \beta>\gamma, \,
                              \beta_{L} \geq \gamma_{L}, \, \alpha \geq \frac{1}{2} \, \max \{\gamma, \gamma_L\} , \\
                              c_4 \, \varepsilon^{-2} \left( \log(\varepsilon) \right)^2 &
                              \text{ if } \beta=\gamma, \, \beta_{L} \geq \gamma_{L}, \, \alpha \geq \frac{1}{2} \, \max \{\gamma, \gamma_L\} , \\
                              c_4 \, \varepsilon^{-2-\frac{\max \{\gamma-\beta,\gamma_{L}-\beta_{L} \}}{\alpha}} &
                              \text{ if } \beta<\gamma, \,
                              \alpha \geq \frac{\max\{\gamma, \gamma_L\}
                               - \max \{ \gamma-\beta, \gamma_L-\beta_L\} }{2},
                             \end{cases}
    \end{equation}
    for some positive constant $c_4$.
\end{The}
%
In order to apply Theorem \ref{Main-Theorem-Giles} and the
multi-level Monte Carlo method, one has to determine the values
$\alpha, \beta, \beta_L >0$ as well as $\gamma, \gamma_L \geq 1$.
Firstly, $\alpha$ denotes the weak order of convergence for the bias
of the finite dimensional approximation $Y^L \in \mathbb{X}_{h_L}$
as the dimension of the approximation subspace increases. This value
is well known for commonly applied approximations $Y^L$. Because the
approximations $(Y^l)_{l \geq 0}$ converge to $X$ in the weak sense,
the differences of two successive approximations $\big(
f(Y^l)-f(Y^{l-1}) \big)_{l \geq 1}$ converge to zero as the
dimensions of the subspaces increase. Then, usually their variances
will also tend to zero with some order $\beta$ and $\beta_L$ for the
approximations applied on levels $0, 1, \ldots, L-1$ and on level
$L$, respectively. Here, we want to point out that
estimates of type 1)--4) in Theorem \ref{Main-Theorem-Giles} are
rather natural and turn out to be no considerable restriction for
typical applications. Finally, the computational costs to evaluate two
correlated approximations $Y^l$ and $Y^{l-1}$ on the finite
dimensional subspaces $\mathbb{X}_{h_l}$ and $\mathbb{X}_{h_{l-1}}$
depend on the dimensions of the subspaces that are proportional to
$h_l^{-1}$. For commonly used discrete time approximations, one
typically has $\gamma=\gamma_L=1$.

%
The calculations for the proof follow the lines of the original
proof due to Giles \cite{Gil08a}. Considering the mean square-error
\begin{equation}
    e(\hat{Y}_{ML}) =
    \left( |\Erw_{\Prob}(f(X))-\Erw_{\tilde{\Prob}}(f(Y^{L}))|^2
    + \Var_{\tilde{\Prob}}(\hat{Y}_{ML}) \right)^{1/2} < \varepsilon
\end{equation}
we make use of the weight $q \in \,]0,1[\,$ and claim that
\begin{equation}
    |\Erw_{\Prob}(f(X))-\Erw_{\tilde{\Prob}}(f(Y^{L}))|^2 < q \, \varepsilon^2 \quad \text{ and } \quad
    \Var_{\tilde{\Prob}}(\hat{Y}_{ML}) < (1-q) \, \varepsilon^2 \, .
\end{equation}
Then, we can calculate $L$ from the bias and we have to solve the
minimization problem
\begin{equation}
    \min_{N_l:0 \leq l \leq L} C(\hat{Y}_{ML})
\end{equation}
under the constraint that $\Var_{\tilde{\Prob}}(\hat{Y}_{ML}) <
(1-q) \, \varepsilon^2$. As a result
of this, we obtain the following
values for $L$ and $N_l$:
\begin{equation}
   L = \left\lceil \frac{\log(q^{-\frac{1}{2}} \, c_{1,\alpha} \, \varepsilon^{-1} \, T^{\alpha})}
   {\alpha \, \log(M)} \right\rceil
\end{equation}
and
$N_0 = \left\lceil \frac{1}{1-q} \, \varepsilon^{-2} \, h_0^{\frac{\beta+\gamma}{2}}
\left( \frac{c_{2,0}}{c_{3,0}} \right)^{\frac{1}{2}} \, \kappa \right\rceil$,
\begin{equation} \label{N_l_Def_Allgemein}
   N_l = \left\lceil \frac{1}{1-q} \, \varepsilon^{-2} \, h_l^{\frac{\beta+\gamma}{2}}
   \left( \frac{c_{2}}{c_{3}} \right)^{\frac{1}{2}} \, \kappa \right\rceil
\end{equation}
for $l=1, \ldots, L-1$ and
$   N_L = \left\lceil \frac{1}{1-q} \, \varepsilon^{-2} \, h_L^{\frac{\beta_{L}+\gamma_{L}}{2}}
   \left( \frac{c_{2,L}}{c_{3,L}} \right)^{\frac{1}{2}} \, \kappa \right\rceil $
for some $q \in \, ]0,1[\,$ where
\begin{itemize}
   \item In case of $\beta>\gamma$ and $\beta_{L} \geq \gamma_{L}$ or
     in case of $\beta<\gamma$ and $\gamma_{L}-\beta_{L} \leq \gamma-\beta$:
     \begin{equation}
         \kappa = \left(c_{2,0} c_{3,0} \right)^{\frac{1}{2}} T^{\frac{\beta-\gamma}{2}}
         + \left( c_2 c_3 \right)^{\frac{1}{2}} \frac{(M^{-1} T)^{\frac{\beta-\gamma}{2}}
         - h_L^{\frac{\beta-\gamma}{2}}}
         {1-M^{\frac{\gamma-\beta}{2}}} + \left( c_{2,L} c_{3,L} \right)^{\frac{1}{2}}
         h_L^{\frac{\beta_{L}-\gamma_{L}}{2}} \, .
     \end{equation}
   \item In case of $\beta=\gamma$ and $\beta_{L} \geq \gamma_{L}$:
     \begin{equation}
         \kappa = \left(c_{2,0} c_{3,0} \right)^{\frac{1}{2}}
         + (L-1) \left( c_{2} c_{3} \right)^{\frac{1}{2}}
         + \left( c_{2,L} c_{3,L} \right)^{\frac{1}{2}} \,
         h_L^{\frac{\beta_{L}-\gamma_{L}}{2}} \, .
     \end{equation}
\end{itemize}
%
%
\section{The improved multi-level Monte Carlo estimator}
\label{Sec3:Improved-MLMC-Estimator}
%
The order of convergence of the multi-level Monte Carlo estimator
$\hat{Y}_{ML}$ given in (\ref{MLMC-estimator-Giles}) is optimal in
the given framework. However, the computational costs can be reduced
if a modified estimator is applied. As yet, the estimator
$\hat{Y}_{ML}$ is based on some weak order $\alpha$ approximations
$Y^l$ for $l=0,1, \ldots, L$ on each level. Now, let us apply some
cheap low order weak approximation $Y^l$ on levels $l=0,1, \ldots,
L-1$ combined with some probably expansive high order weak
approximation $\check{Y}^L$ on the finest level $L$. The idea is,
that the approximations $Y^l$ contribute
a variance reduction
while the approximation $\check{Y}^L$ results in a small bias of the
multi-level Monte Carlo estimator, thus reducing the number of
levels needed to attain a prescribed accuracy.

Let $Y$ be an order $\alpha$ weak approximation method and let
$\check{Y}$ be an order $p$ weak approximation method applied on the
finest level. Further, let $L=L_p$ with
\begin{equation}
   L_p = \left\lceil \frac{\log(q^{-\frac{1}{2}} \, c_{1,p} \, \varepsilon^{-1} \, T^{p})}
   {p \, \log(M)} \right\rceil
\end{equation}
denote the number of levels in order to indicate the dependence on
the weak order $p$. Then, we define the modified multi-level Monte
Carlo estimator
by
\begin{equation} \label{Sec2:modified-MLMC-estimator}
    \hat{Y}_{ML(\alpha,p)} = \sum_{l=0}^{L_p} \hat{Y}^l
\end{equation}
with the estimators $\hat{Y}^l$ for $l=0,1, \ldots, L_p-1$ based on
the order $\alpha$ weak approximations $Y^l$ as defined in
Section~\ref{Section2:MLMC-Simulation-Original}, however now
applying the modified estimator
\begin{equation}
    \hat{Y}^{L_p} = \frac{1}{N_{L_p}} \sum_{i=1}^{N_{L_p}}
    \left( f({\check{Y}^{L_p}})^{(i)} - f({Y^{L_p-1}})^{(i)} \right)
\end{equation}
which combines the weak order $\alpha$ approximations $Y^{L_p-1}$
with the weak order $p$ approximations $\check{Y}^{L_p}$. Clearly,
all conditions of Theorem~\ref{Main-Theorem-Giles} have to be
fulfilled for $Y^L$ replaced by $\check{Y}^L$. Then, in the case of
$p > \alpha$, the improved multi-level Monte Carlo estimator
$\hat{Y}_{ML(\alpha,p)}$ features significantly reduced
computational costs compared to the originally proposed estimator
$\hat{Y}_{ML}=\hat{Y}_{ML(\alpha,\alpha)}$.
\begin{Sat} \label{Main-Prop-Improvement}
    Let conditions 1)--7) of Theorem~\ref{Main-Theorem-Giles}
    be fulfilled and suppose that there exist constants $\hat{c}_{3,0}, \hat{c}_3,
    \hat{c}_{3,L_p}, \delta_i > 0$ and $\hat{c}_{3,0}^{(i)}, \hat{c}_3^{(i)},
    \hat{c}_{3,L_p}^{(i)} \geq 0$
    such that for the computational costs
    \begin{enumerate}
        \item[5')] $C(Y^0) = \hat{c}_{3,0} \, T \, h_0^{-\gamma}
        + \sum_{i=1}^{k} \hat{c}_{3,0}^{(i)} \, T \, h_0^{-\gamma+\delta_i}$,
        \item[6')] $C(Y^l,Y^{l-1}) = \hat{c}_{3} \, T \, h_l^{-\gamma}
        + \sum_{i=1}^{k} \hat{c}_{3}^{(i)} \, T \, h_l^{-\gamma+\delta_i}$ \quad for $l=1, \ldots, L_p-1$,
        \item[7')] $C(\check{Y}^{L_p},Y^{L_p-1}) = \hat{c}_{3,L_p} \, T \, h_{L_p}^{-\gamma_{L_p}}
        + \sum_{i=1}^{k} \hat{c}_{3,L_p}^{(i)} \, T \, h_{L_p}^{-\gamma_{L_p}+\delta_i}$
    \end{enumerate}
    with some $\gamma, \gamma_{L_p} \geq 1$ such that
    $\gamma-\delta_i \geq 1$ and $\gamma_{L_p}-\delta_i \geq 1$.
    Then, the multi-level Monte Carlo estimator $\hat{Y}_{ML(\alpha,p)}$ based on
    a weak order $\alpha>0$ approximation scheme on levels $0,1, \ldots, L_p-1$ and
    some weak order $p > \alpha$ approximation scheme on level $L_p$ has reduced computational
    costs:
    \begin{enumerate}[i)]
       \item In case of $\beta>\gamma$ and $\beta-\gamma<\beta_{L_p}-\gamma_{L_p}$, there exists some
             $\varepsilon_0>0$ such that for all $\varepsilon \in \, ]0,\varepsilon_0]$ it holds
             \begin{equation} \label{Main-Prop-Improvement-Aussage1}
                \frac{C(\hat{Y}_{ML(\alpha,\alpha)})(\varepsilon)}
                {C(\hat{Y}_{ML(\alpha,p)})(\varepsilon)} > 1
             \end{equation}
             provided that $\alpha \geq \frac{\gamma}{2}$,
             $p \geq \frac{1}{2} \max \{\gamma, \gamma_{L_p}\}$ and
             $p>\tfrac{1}{4} \max \{\beta+\gamma,\beta-\gamma+2\gamma_{L_p} \}$.
             In case of $\beta>\gamma$ and $\beta-\gamma=\beta_{L_p}-\gamma_{L_p}$ then
             (\ref{Main-Prop-Improvement-Aussage1}) holds
             if in addition $c_2 c_3 > (1-M^{\frac{\gamma-\beta}{2}})^2 c_{2,L_p} c_{3,L_p}$ and
             $
             \hat{c}_3^2
             \frac{c_2}{c_3} > (1-M^{\frac{\gamma-\beta}{2}})^2
             \hat{c}_{3,L_p}^2
             \frac{c_{2,L_p}}{c_{3,L_p}}$.
             Further, for $0<\beta-\gamma \leq \beta_{L_p}-\gamma_{L_p}$
             it holds $C(\hat{Y}_{ML(\alpha,p)})(\varepsilon)
             = O( \varepsilon^{-2} )$ if $\alpha>0$ and $p \geq \tfrac{1}{2} \max \{\gamma,
             \gamma_{L_p}\}$.
       \item In case of $\beta=\gamma$ and $\beta_{L_p} \geq \gamma_{L_p}$ and if
             $p \geq \frac{1}{2} \max \{\gamma, \gamma_{L_p} \}$, $\alpha \geq \frac{\gamma}{2}$, it holds
             \begin{equation} \label{Main-Prop-Improvement-Aussage2}
                \lim_{\varepsilon \to 0} \frac{C(\hat{Y}_{ML(\alpha,\alpha)})(\varepsilon)}
                {C(\hat{Y}_{ML(\alpha,p)})(\varepsilon)} \geq \left( \frac{p}{\alpha}
                \right)^2
             \end{equation}
             and $C(\hat{Y}_{ML(\alpha,p)})(\varepsilon)
             = O( \varepsilon^{-2} ( \log(\varepsilon) )^2 )$ if $\alpha>0$
             and $p \geq \frac{1}{2} \max \{\gamma, \gamma_{L_p}\}$.
       \item In case of $\beta<\gamma$ and $\gamma-\beta=\gamma_{L_p}-\beta_{L_p}$ it holds
             \begin{equation} \label{Main-Prop-Improvement-Aussage3}
                \begin{split}
                \lim_{\varepsilon \to 0} \frac{C(\hat{Y}_{ML(p,p)})(\varepsilon)}
                {C(\hat{Y}_{ML(\alpha,p)})(\varepsilon)}
                \geq & \,
                M^{2(\gamma-\beta)} \left( \frac{\hat{c}_3
                c_2}{\hat{c}_{3,L_p}
                c_{2,L_p}} + \frac{\hat{c}_3
                (c_2 c_{3,L_p})^{1/2}}{\hat{c}_{3,L_p}
                (c_{2,L_p} c_3)^{1/2}}
                \left( M^{\frac{\gamma-\beta}{2}} -1 \right) \right. \\
                & \left. + \left( \frac{c_2 c_3}{c_{2,L_p} c_{3,L_p}} \right)^{1/2}
                \left( M^{\frac{\gamma-\beta}{2}} -1 \right) + \left( M^{\frac{\gamma-\beta}{2}} -1 \right)^2 \right)^{-1}
                \end{split}
             \end{equation}
             if $p > \frac{1}{2} (\max\{\gamma, \gamma_{L_p}\} - \gamma+\beta)$.
             If the pa\-ra\-me\-ter $q \in \,]0,1[\,$
             is chosen as
             \begin{equation} \label{Main-Prop-Improvement-Aussage4}
                q = \frac{\gamma-\beta}{\gamma-\beta+2p}
             \end{equation}
             then the computational costs
             $C(\hat{Y}_{ML(\alpha,p)})$ are asymptotically minimal.
             In general, if $\beta<\gamma$ or if $\beta_{L_p}<\gamma_{L_p}$ then it holds
             $C(\hat{Y}_{ML(\alpha,p)})(\varepsilon)
             = O \big( \varepsilon^{-2-\frac{\max \{\gamma-\beta,\gamma_{L_p}-\beta_{L_p} \}}{p}} \big)$
             for $p \geq \frac{1}{2} (\max\{\gamma, \gamma_{L_p}\} -
             \min \{\gamma-\beta,\gamma_{L_p}-\beta_{L_p}\})$.
    \end{enumerate}
\end{Sat}
%
We note, that in relations 5')--7') of
Proposition~\ref{Main-Prop-Improvement} a more detailed polynomial
dependence of the computational costs from the dimension of the
approximation subspaces has to be taken into account. E.g., standard
discrete time approximation methods
possess polynomial computational costs and the constants are known
explicitly.
%
\begin{proof}
In the following, we will first state some basic formulas and conditions used
in the remaining part of
the proof. Then we will calculate lower and upper bounds for
the computational costs in the case $\beta\neq\gamma$. Those will then be used to prove
first i) and then iii). Finally, case ii) with $\beta=\gamma$ is considered.

\emph{Basic formulas.}
Assume that $\varepsilon<1$.
Let $\delta_0=0$, $\hat{c}_{3,0}^{(0)}=\hat{c}_{3,0}$, $\hat{c}_{3}^{(0)}=\hat{c}_{3}$
and $\hat{c}_{3,L_p}^{(0)}=\hat{c}_{3,L_p}$. Then, the computational costs for $\hat{Y}_{ML(\alpha,p)}$ are
\begin{equation}
    \begin{split}
    C(\hat{Y}_{ML(\alpha,p)}) = & \, \sum_{i=0}^k \hat{c}_{3,0}^{(i)} \, T \, h_0^{-\gamma+\delta_i}
    \, N_0 + \sum_{i=0}^k \sum_{l=1}^{L_p-1} \hat{c}_3^{(i)} \, T \, h_l^{-\gamma+\delta_i}
    \, N_l \\
    & \, + \sum_{i=0}^k \hat{c}_{3,L_p}^{(i)} \, T \, h_{L_p}^{-\gamma_{L_p}+\delta_i}
    \, N_{L_p}
    \end{split}
\end{equation}
with $L=L_p = \left\lceil \frac{\log(q^{-\frac{1}{2}} \, c_{1,p} \, \varepsilon^{-1} \, T^p)}
{p \, \log(M)} \right\rceil$ and $N_l$ for $l=0,1, \ldots, L_p$ given in \eqref{N_l_Def_Allgemein}.
Without loss of generality, suppose that $\delta_i \neq \delta_j$
for $i \neq j$ and that $\delta_{k}=\tfrac{\gamma-\beta}{2}$ with
$\hat{c}_{3,0}^{(k)} = \hat{c}_3^{(k)} = \hat{c}_{3,L_p}^{(k)} = 0$
in the case of $\beta \geq \gamma$.
In the following, we make use of the two estimates
\begin{align}
    L_{\alpha} &\geq \frac{\log(\varepsilon^{-1})}{\alpha \, \log(M)}
    + \frac{\log(q^{-\frac{1}{2}} \, c_{1,\alpha} \,
    T^{\alpha})}{\alpha \, \log(M)} \, ,
    \label{L-lower-bound} \\
    L_p - 1 &\leq \frac{\log(\varepsilon^{-1})}{p \log(M)}
    + \frac{\log(q^{-\frac{1}{2}} c_{1,p} T^p)}{p \log(M)} \, .
    \label{L-1-upper-bound}
\end{align}

\emph{Lower bound for $\beta\neq\gamma$.}
Let $\beta \neq \gamma$. Then, we obtain the lower bound
{\allowdisplaybreaks
\begin{align}
    C(\hat{Y}_{ML(\alpha,\alpha)})(\varepsilon)
    \geq &\, \frac{T \kappa \, \varepsilon^{-2}}{1-q}
    \sum_{i=0}^k \left(
    h_0^{\frac{\beta-\gamma}{2}+\delta_i} \hat{c}_{3,0}^{(i)}
    \left(\frac{c_{2,0}}{c_{3,0}} \right)^{1/2}
    + \sum_{l=1}^{L_{\alpha}} h_l^{\frac{\beta-\gamma}{2}+\delta_i}
    \hat{c}_3^{(i)} \left( \frac{c_{2}}{c_3} \right)^{1/2}
    \right) \notag \\
    \geq& \, \frac{T}{1-q} \, \varepsilon^{-2} \Bigg[
    \sum_{i=0}^k T^{\beta-\gamma+\delta_i} \hat{c}_{3,0}^{(i)} c_{2,0}\notag \\
    & + \sum_{i=0}^k T^{\frac{\beta-\gamma}{2}+\delta_i} \hat{c}_{3,0}^{(i)}
    \left( \frac{c_{2,0} c_2 c_3}{c_{3,0}} \right)^{1/2}
    \frac{ T^{\frac{\beta-\gamma}{2}} - h_{L_{\alpha}}^{\frac{\beta-\gamma}{2}}}{M^{\frac{\beta-\gamma}{2}}-1} \notag \\
    & + \sum_{i=0}^{k-1} \hat{c}_{3}^{(i)} \left( \frac{c_2 c_{2,0} c_{3,0}}{c_3} \right)^{1/2} T^{\frac{\beta-\gamma}{2}} \cdot
    \frac{T^{\frac{\beta-\gamma}{2}+\delta_i}
    - h_{L_{\alpha}}^{\frac{\beta-\gamma}{2}+\delta_i}}{M^{\frac{\beta-\gamma}{2}+\delta_i}-1} \notag \\
    & + \hat{c}_3^{(k)} \left( \frac{c_2 c_{2,0} c_{3,0}}{c_3} \right)^{1/2}
    T^{\frac{\beta-\gamma}{2}} \left( \frac{\log(\varepsilon^{-1})}{\alpha \log(M)} + \frac{\log(q^{-\frac{1}{2}}
    c_{1,\alpha} T^{\alpha})}{\alpha \log(M)} \right) \notag \\
    & + \sum_{i=0}^{k-1} \hat{c}_{3}^{(i)} c_2 \frac{T^{\frac{\beta-\gamma}{2}+\delta_i}
    - h_{L_{\alpha}}^{\frac{\beta-\gamma}{2}+\delta_i}}{M^{\frac{\beta-\gamma}{2}+\delta_i}-1} \cdot
    \frac{T^{\frac{\beta-\gamma}{2}} - h_{L_{\alpha}}^{\frac{\beta-\gamma}{2}}}{M^{\frac{\beta-\gamma}{2}}-1} \notag \\
    &+ \hat{c}_3^{(k)} c_2
    \frac{T^{\frac{\beta-\gamma}{2}} - h_{L_{\alpha}}^{\frac{\beta-\gamma}{2}}}{M^{\frac{\beta-\gamma}{2}}-1}
    \left( \frac{\log(\varepsilon^{-1})}{\alpha \log(M)} + \frac{\log(q^{-\frac{1}{2}}
    c_{1,\alpha} T^{\alpha})}{\alpha \log(M)} \right)
    \Bigg] \label{Proof-Main-Prop-Lower-Bound-Ieq-alg}
\end{align}
}
where $\hat{c}_{3,L_{\alpha}}^{(i)} = \hat{c}_{3}^{(i)}$,
$c_{2,L_{\alpha}}=c_2$, $c_{3,L_{\alpha}}=c_3$, $\beta_{L_{\alpha}}=\beta$
and $\gamma_{L_{\alpha}}=\gamma$ for $\hat{Y}_{ML(\alpha,\alpha)}$.

\emph{Upper bound for $\beta\neq\gamma$.}
Next, we calculate for the case of $\beta \neq \gamma$ the upper
bound
{\allowdisplaybreaks
\begin{align}
    C(\hat{Y}_{ML(\alpha,p)})(\varepsilon)\notag\\ \leq &\, \frac{T \kappa \, \varepsilon^{-2}}{1-q}
    \sum_{i=0}^k \left(
    h_0^{\frac{\beta-\gamma}{2}+\delta_i} \hat{c}_{3,0}^{(i)} \left(\frac{c_{2,0}}{c_{3,0}} \right)^{1/2}
    + \sum_{l=1}^{L_p-1} h_l^{\frac{\beta-\gamma}{2}+\delta_i}
    \hat{c}_3^{(i)} \left( \frac{c_{2}}{c_3} \right)^{1/2} \right. \notag \\
    & \left. \, + h_{L_p}^{\frac{\beta_{L_p}-\gamma_{L_p}}{2}+\delta_i}
    \hat{c}_{3,L_p}^{(i)} \left( \frac{c_{2,L_p}}{c_{3,L_p}} \right)^{1/2}
    \right) \notag \\
    & + T \sum_{i=0}^k \left( \hat{c}_{3,0}^{(i)} h_0^{-\gamma+\delta_i}
    + \hat{c}_3^{(i)} \sum_{l=1}^{L_p-1} h_l^{-\gamma+\delta_i}
    + \hat{c}_{3,L_p}^{(i)} h_{L_p}^{-\gamma_{L_p}+\delta_i} \right) \notag \\
    \leq & \, \frac{T}{1-q} \, \varepsilon^{-2} \Bigg[
    \sum_{i=0}^k \hat{c}_{3,0}^{(i)} \left( c_{2,0} T^{\beta-\gamma+\delta_i}
    +
    \left( \frac{c_{2,0} c_2 c_3}{c_{3,0}}
    \right)^{1/2} \Lambda_0 T^{\frac{\beta-\gamma}{2}+\delta_i} \right. \notag \\
    & + \left.
    \left( \frac{c_{2,0} c_{2,L_p} c_{3,L_p}}{c_{3,0}}
    \right)^{1/2} T^{\frac{\beta-\gamma}{2}+\delta_i} h_{L_p}^{\frac{\beta_{L_p}-\gamma_{L_p}}{2}}
    \right) \notag \\
    & + \sum_{i=0}^{k-1} \hat{c}_3^{(i)} \left( \frac{c_{2,0} c_{3,0} c_2}{c_3} \right)^{1/2}
    T^{\frac{\beta-\gamma}{2}} \Lambda_i
    + \sum_{i=0}^{k-1} \hat{c}_3^{(i)} c_2 \Lambda_i \Lambda_0 \notag \\
    & + \left( \hat{c}_3^{(k)} c_2 \Lambda_0
    + \hat{c}_3^{(k)} \left( \frac{c_{2} c_{2,L_p} c_{3,L_p}}{c_3} \right)^{1/2}
    h_{L_p}^{\frac{\beta_{L_p}-\gamma_{L_p}}{2}} \right. \notag \\
    & + \left. \hat{c}_3^{(k)} \left( \frac{c_{2,0} c_{3,0} c_2}{c_3} \right)^{1/2}
    T^{\frac{\beta-\gamma}{2}} \right) \left(
    \frac{\log(\varepsilon^{-1})}{p \log(M)}
    + \frac{\log(q^{-\frac{1}{2}} c_{1,p} T^p)}{p \log(M)} \right) \notag \\
    & + \sum_{i=0}^{k-1} \hat{c}_3^{(i)} \left( \frac{c_{2} c_{2,L_p} c_{3,L_p}}{c_3} \right)^{1/2}
    \Lambda_i h_{L_p}^{\frac{\beta_{L_p}-\gamma_{L_p}}{2}} \notag \\
    & + \sum_{i=0}^k \hat{c}_{3,L_p}^{(i)}
    \left( \frac{c_{2,0} c_{3,0} c_{2,L_p}}{c_{3,L_p}} \right)^{1/2} T^{\frac{\beta-\gamma}{2}}
    h_{L_p}^{\frac{\beta_{L_p}-\gamma_{L_p}}{2}+\delta_i} \notag \\
    &+ \sum_{i=0}^k \hat{c}_{3,L_p}^{(i)} \left(
    \left( \frac{c_{2} c_{3} c_{2,L_p}}{c_{3,L_p}} \right)^{1/2} \Lambda_0
    h_{L_p}^{\frac{\beta_{L_p}-\gamma_{L_p}}{2}+\delta_i}
    +
    c_{2,L_p} h_{L_p}^{\beta_{L_p}-\gamma_{L_p}+\delta_i} \right)  \Bigg] \notag \\
    & + T \sum_{i=0}^k \left( \hat{c}_{3,0}^{(i)} T^{\delta_i-\gamma}
    + \hat{c}_3^{(i)} \frac{(M^{-1} T)^{\delta_i-\gamma} - h_{L_p}^{\delta_i-\gamma}}
    {1-M^{\gamma-\delta_i}} + \hat{c}_{3,L_p}^{(i)} h_{L_p}^{\delta_i-\gamma_{L_p}} \right)
    \label{Proof-Main-Prop-Upper-Bound-Ieq-alg}
\end{align}
}
with $\Lambda_i = \frac{(M^{-1}
T)^{\frac{\beta-\gamma}{2}+\delta_i} -
h_{L_p}^{\frac{\beta-\gamma}{2}+\delta_i}}{1-M^{\frac{\gamma-\beta}{2}-\delta_i}}$
for $i=0,\dots,k-1$.

\emph{Proof of i).}
In case of $\beta>\gamma$ and $\beta_{L_p} > \gamma_{L_p}$, we
prove that there exists some $\varepsilon_0>0$ such that for all
$\varepsilon \in \, ]0,\varepsilon_0]$ it follows
$C(\hat{Y}_{ML(\alpha,\alpha)})(\varepsilon) >
C(\hat{Y}_{ML(\alpha,p)})(\varepsilon)$.
From the lower bound (\ref{Proof-Main-Prop-Lower-Bound-Ieq-alg}) for
$C(\hat{Y}_{ML(\alpha,\alpha)})(\varepsilon)$ and the upper bound
(\ref{Proof-Main-Prop-Upper-Bound-Ieq-alg}) for
$C(\hat{Y}_{ML(\alpha,p)})(\varepsilon)$ we get the estimate
{\allowdisplaybreaks
\begin{align}
    & C(\hat{Y}_{ML(\alpha,\alpha)})(\varepsilon) - C(\hat{Y}_{ML(\alpha,p)})(\varepsilon) \notag \\
    & \geq \frac{T}{1-q} \, \varepsilon^{-2} 
    \left( \sum_{i=0}^{k-1} T^{\frac{\beta-\gamma}{2}+\delta_i} \hat{c}_{3,0}^{(i)}
    \left( \frac{c_{2,0} c_2 c_3}{c_{3,0}} \right)^{1/2}
    \frac{h_{L_p}^{\frac{\beta-\gamma}{2}} - M^{\frac{\gamma-\beta}{2}}
    h_{L_{\alpha}}^{\frac{\beta-\gamma}{2}}}{1-M^{\frac{\gamma-\beta}{2}}} \right. \notag \\
    & + \sum_{i=0}^{k-1} \hat{c}_{3}^{(i)} \left( \frac{c_2 c_{2,0} c_{3,0}}{c_3} \right)^{1/2} T^{\frac{\beta-\gamma}{2}} \cdot
    \frac{h_{L_p}^{\frac{\beta-\gamma}{2}+\delta_i}
    - M^{\frac{\gamma-\beta}{2}-\delta_i} h_{L_{\alpha}}^{\frac{\beta-\gamma}{2}+\delta_i}}{1-M^{\frac{\gamma-\beta}{2}-\delta_i}} \notag \\
    & + \sum_{i=0}^{k-1} \hat{c}_{3}^{(i)} c_2 \left(
    \frac{(M^{-1} T)^{\frac{\beta-\gamma}{2}+\delta_i} \left( h_{L_p}^{\frac{\beta-\gamma}{2}}
    - M^{\frac{\gamma-\beta}{2}} h_{L_{\alpha}}^{\frac{\beta-\gamma}{2}} \right)
     }
    {(1-M^{\frac{\gamma-\beta}{2}-\delta_i}) (1-M^{\frac{\gamma-\beta}{2}}) } \right. \notag \\
    & + \left. \frac{(M^{-1} T)^{\frac{\beta-\gamma}{2}} \left( h_{L_p}^{\frac{\beta-\gamma}{2}+\delta_i}
    - M^{\frac{\gamma-\beta}{2}-\delta_i} h_{L_{\alpha}}^{\frac{\beta-\gamma}{2}+\delta_i} \right)
    - h_{L_p}^{\beta-\gamma+\delta_i} + M^{\gamma-\beta-\delta_i} h_{L_{\alpha}}^{\beta-\gamma+\delta_i} }
    {(1-M^{\frac{\gamma-\beta}{2}-\delta_i}) (1-M^{\frac{\gamma-\beta}{2}}) } \right) \notag \\
    & - \sum_{i=0}^{k-1} \hat{c}_{3,0}^{(i)}
    \left( \frac{c_{2,0} c_{2,L_p} c_{3,L_p}}{c_{3,0}}
    \right)^{1/2} T^{\frac{\beta-\gamma}{2}+\delta_i} h_{L_p}^{\frac{\beta_{L_p}-\gamma_{L_p}}{2}}
    \notag \\
    & - \sum_{i=0}^{k-1} \hat{c}_3^{(i)} \left( \frac{c_{2} c_{2,L_p} c_{3,L_p}}{c_3} \right)^{1/2}
    \Lambda_i h_{L_p}^{\frac{\beta_{L_p}-\gamma_{L_p}}{2}} \notag \\
    & - \sum_{i=0}^{k-1} \hat{c}_{3,L_p}^{(i)}
    \left( \frac{c_{2,0} c_{3,0} c_{2,L_p}}{c_{3,L_p}} \right)^{1/2} T^{\frac{\beta-\gamma}{2}}
    h_{L_p}^{\frac{\beta_{L_p}-\gamma_{L_p}}{2}+\delta_i} \notag \\
    & \left. - \sum_{i=0}^{k-1} \hat{c}_{3,L_p}^{(i)} \left(
    \left( \frac{c_{2} c_{3} c_{2,L_p}}{c_{3,L_p}} \right)^{1/2} \Lambda_0
    h_{L_p}^{\frac{\beta_{L_p}-\gamma_{L_p}}{2}+\delta_i}
    + 
    c_{2,L_p} h_{L_p}^{\beta_{L_p}-\gamma_{L_p}+\delta_i} \right) \right) \notag \\
    & - T \sum_{i=0}^{k-1} \left( \hat{c}_{3,0}^{(i)} T^{\delta_i-\gamma}
    + \hat{c}_3^{(i)} \frac{(M^{-1} T)^{\delta_i-\gamma} - h_{L_p}^{\delta_i-\gamma}}
    {1-M^{\gamma-\delta_i}} + \hat{c}_{3,L_p}^{(i)} h_{L_p}^{\delta_i-\gamma_{L_p}} \right) .
    \label{Proof-Main-Prop-Difference-Ieq-alg}
\end{align}
}
\noindent
In the following, we make use of the estimates $M^{-1}
c_{1,\alpha}^{-\frac{1}{\alpha}} q^{\frac{1}{2 \alpha}}
\varepsilon^{\frac{1}{\alpha}} \leq h_{L_{\alpha}} \leq
c_{1,\alpha}^{-\frac{1}{\alpha}} q^{\frac{1}{2 \alpha}}
\varepsilon^{\frac{1}{\alpha}}$ and $M^{-1} c_{1,p}^{-\frac{1}{p}}
q^{\frac{1}{2 p}} \varepsilon^{\frac{1}{p}} \leq h_{L_{p}} \leq
c_{1,p}^{-\frac{1}{p}} q^{\frac{1}{2 p}} \varepsilon^{\frac{1}{p}}$,
i.e.\ we have $h_{L_{p}} \to 0$ and $h_{L_{\alpha}} \to 0$ as
$\varepsilon \to 0$.

Multiplying both sides of (\ref{Proof-Main-Prop-Difference-Ieq-alg})
with $\frac{1-q}{T} \, \varepsilon^2 \, h_{L_p}^{-\frac{\min
\{\beta-\gamma,\beta_{L_p}-\gamma_{L_p}\}}{2}}$ and taking into
account the assumptions $4p>\beta+\gamma$ and
$4p>\beta-\gamma+2\gamma_{L_p}$ results in {\allowdisplaybreaks
\begin{align}
    & \frac{1-q}{T} \, \varepsilon^2 \,
    h_{L_p}^{-\frac{\min \{\beta-\gamma,\beta_{L_p}-\gamma_{L_p}\}}{2}}
    \left( C(\hat{Y}_{ML(\alpha,\alpha)})(\varepsilon) - C(\hat{Y}_{ML(\alpha,p)})(\varepsilon)
    \right) \notag \\
    & \geq \Bigg[
    \sum_{i=0}^{k-1} T^{\frac{\beta-\gamma}{2}+\delta_i} \hat{c}_{3,0}^{(i)}
    \left( \frac{c_{2,0}}{c_{3,0}} \right)^{1/2} \left( \left( c_2 c_3 \right)^{1/2}
    \frac{h_{L_p}^{\frac{\beta-\gamma}{2}} }{1-M^{\frac{\gamma-\beta}{2}}}
    - \left( c_{2,L_p} c_{3,L_p} \right)^{1/2} h_{L_p}^{\frac{\beta_{L_p}-\gamma_{L_p}}{2}}
    \right) \notag \\
    & + T^{\frac{\beta-\gamma}{2}} (c_{2,0} c_{3,0})^{1/2}
    \left( \hat{c}_{3}^{(0)} \left( \frac{c_2}{c_3} \right)^{1/2}
    \frac{h_{L_p}^{\frac{\beta-\gamma}{2}}}{1-M^{\frac{\gamma-\beta}{2}}}
    - \hat{c}_{3,L_p}^{(0)} \left( \frac{c_{2,L_p}}{c_{3,L_p}} \right)^{1/2} h_{L_p}^{\frac{\beta_{L_p}-\gamma_{L_p}}{2}} \right) \notag \\
    & + \sum_{i=0}^{k-1} \hat{c}_{3}^{(i)} {c_2}^{1/2}
    \frac{(M^{-1} T)^{\frac{\beta-\gamma}{2}+\delta_i} }
    {1-M^{\frac{\gamma-\beta}{2}-\delta_i} } \left( (c_2)^{1/2}
    \frac{h_{L_p}^{\frac{\beta-\gamma}{2}} }{1-M^{\frac{\gamma-\beta}{2}} }
    - \left( \frac{c_{2,L_p} c_{3,L_p}}{c_3} \right)^{1/2}
    h_{L_p}^{\frac{\beta_{L_p}-\gamma_{L_p}}{2}} \right)
    \notag \\
    & + \frac{(M^{-1} T)^{\frac{\beta-\gamma}{2}} }
    {1-M^{\frac{\gamma-\beta}{2}} } {c_2}^{1/2} \left( \hat{c}_{3}^{(0)} {c_2}^{1/2}
    \frac{h_{L_p}^{\frac{\beta-\gamma}{2}} }
    {1-M^{\frac{\gamma-\beta}{2}} }
    - \hat{c}_{3,L_p}^{(0)} \left( \frac{c_{3} c_{2,L_p}}{c_{3,L_p}} \right)^{1/2}
    h_{L_p}^{\frac{\beta_{L_p}-\gamma_{L_p}}{2}} \right) \notag \\
    &
    +
    o \left( h_{L_p}^{\frac{\min \{\beta_{L_p}-\gamma_{L_p}, \beta-\gamma \}}{2}} \right)
    \Bigg]
    h_{L_p}^{-\frac{\min \{\beta-\gamma,\beta_{L_p}-\gamma_{L_p}\}}{2}}\,.
    \label{Proof-Main-Prop-Difference2-Ieq-alg}
\end{align}
}
As a result of (\ref{Proof-Main-Prop-Difference2-Ieq-alg}) it follows that
in the case of $\beta-\gamma<\beta_{L_p}-\gamma_{L_p}$ there exists some
$\varepsilon_0>0$ such that
\begin{equation} \label{Proof-Main-Prop-Result-case-1}
    \frac{C(\hat{Y}_{ML(\alpha,\alpha)})(\varepsilon)}
    {C(\hat{Y}_{ML(\alpha,p)})(\varepsilon)} > 1
\end{equation}
for all $\varepsilon \in \, ]0,\varepsilon_0]$. In the case of
$\beta-\gamma=\beta_{L_p}-\gamma_{L_p}$ there exists some
$\varepsilon_0>0$ such that (\ref{Proof-Main-Prop-Result-case-1})
holds for all $\varepsilon \in \, ]0,\varepsilon_0]$ if $c_2 c_3 >
(1-M^{\frac{\gamma-\beta}{2}})^2 c_{2,L_p} c_{3,L_p}$ and $\left(
\hat{c}_3^{(0)} \right)^2 \frac{c_2}{c_3} >
(1-M^{\frac{\gamma-\beta}{2}})^2 \left( \hat{c}_{3,L_p}^{(0)}
\right)^2 \frac{c_{2,L_p}}{c_{3,L_p}}$. Finally,
$C(\hat{Y}_{ML(\alpha,p)})(\varepsilon) = O( \varepsilon^{-2} )$
follows from (\ref{Proof-Main-Prop-Upper-Bound-Ieq-alg}).
%
%

\emph{Proof of iii).}
In case of $\beta<\gamma$ and $\beta<2p$, we have to compare the
dominating terms as $\varepsilon \to 0$. Therefore, we get from the
lower bound that
\begin{align}
    C(\hat{Y}_{ML(p,p)})(\varepsilon) \geq \,
    & \frac{q^{\frac{\beta-\gamma}{2p}}}{1-q} \,
    \varepsilon^{-2-\frac{\gamma-\beta}{p}} \, T \,
    \hat{c}_{3,L_p}^{(0)} \, c_{2,L_p} \,
    c_{1,p}^{\frac{\gamma-\beta}{p}}
    \, M^{\gamma-\beta} \left( M^{\frac{\beta-\gamma}{2}}-1
    \right)^{-2} \notag \\
    & + o(\varepsilon^{-2-\frac{\gamma-\beta}{p}}) \label{Proof-Main-Prop-Lower-Bound-Ieq-simp}
\end{align}
and from the upper bound
\begin{align}
    C(\hat{Y}_{ML(\alpha,p)})(\varepsilon) \leq \, &
    \frac{q^{\frac{\beta-\gamma}{2p}}}{1-q} \,
    \varepsilon^{-2-\frac{\gamma-\beta}{p}} \, T \,
    c_{1,p}^{\frac{\gamma-\beta}{p}} \left(
    \frac{\hat{c}_{3}^{(0)} c_{2}}{\left( 1-M^{\frac{\gamma-\beta}{2}} \right)^{2}}
    - \frac{\hat{c}_{3}^{(0)} \left( c_{2} c_{2,L_p} c_{3,L_p} \right)^{1/2}}{c_3^{1/2}
    \left( 1-M^{\frac{\gamma-\beta}{2}} \right)} \right. \notag \\
    & - \left. \frac{\hat{c}_{3,L_p}^{(0)} \left( c_{2} c_3 c_{2,L_p} \right)^{1/2}}{c_{3,L_p}^{1/2}
    \left( 1-M^{\frac{\gamma-\beta}{2}} \right)} + \hat{c}_{3,L_p}^{(0)} c_{2,L_p} \right)
    + o(\varepsilon^{-2-\frac{\gamma-\beta}{p}}) . \label{Proof-Main-Prop-Upper-Bound-Ieq-simp}
\end{align}
Making use of these two estimates
(\ref{Proof-Main-Prop-Lower-Bound-Ieq-simp}) and
(\ref{Proof-Main-Prop-Upper-Bound-Ieq-simp}), this results in the
estimate (\ref{Main-Prop-Improvement-Aussage3})
where $\beta_{L_p}<\gamma_{L_p}$ because we
require that $\beta_{L_p}-\gamma_{L_p}=\beta-\gamma<0$.

In general, it follows that $C(\hat{Y}_{ML(\alpha,p)})(\varepsilon)
= O \big( \varepsilon^{-2-\frac{\max
\{\gamma-\beta,\gamma_{L_p}-\beta_{L_p} \}}{p}} \big)$ due to
the upper bound (\ref{Proof-Main-Prop-Upper-Bound-Ieq-alg}) for
$\beta<\gamma$ and any $\beta_{L_p}>0$, $\gamma_{L_p} \geq 1$.
Further, there is an asymptotically optimal choice for the parameter
$q \in \,]0,1[\,$ such that the computational costs are
asymptotically minimal. Calculating a lower bound for
$C(\hat{Y}_{ML(\alpha,p)})(\varepsilon)$ and taking into account the
upper bound (\ref{Proof-Main-Prop-Upper-Bound-Ieq-simp}), we
get
\begin{equation}
    C(\hat{Y}_{ML(\alpha,p)})(\varepsilon) =
    \frac{1}{1-q} \, \varepsilon^{-2-\frac{\gamma-\beta}{p}} \,
    q^{\frac{\beta-\gamma}{2p}}
    \, C + o(\varepsilon^{-2-\frac{\gamma-\beta}{p}})
\end{equation}
with some constant $C>0$ independent of $q$ and $\varepsilon$.
Now, we have to find some $\hat{q} \in \,]0,1[\,$ such that
\begin{equation}
    C \varepsilon^{-2-\frac{\gamma-\beta}{p}}
    \frac{\hat{q}^{\frac{\beta-\gamma}{2p}}}{1-\hat{q}}
    = \min_{q \in \, ]0,1[ \,} C \varepsilon^{-2-\frac{\gamma-\beta}{p}}
    \frac{q^{\frac{\beta-\gamma}{2p}}}{1-q}
\end{equation}
for all $0<\varepsilon<1$. Solving this minimization problem leads
to
\begin{equation}
    \hat{q} = \frac{\gamma-\beta}{\gamma-\beta+2p}
\end{equation}
which is asymptotically the optimal choice for $q \in \, ]0,1[\,$ in
case of $\beta<\gamma$.
%
%

\emph{Lower bound for $\beta=\gamma$.}
In case of $\beta=\gamma$, we get the following lower bound
{\allowdisplaybreaks
\begin{align}
    C(\hat{Y}_{ML(\alpha,\alpha)})(\varepsilon) \geq &\, \frac{T}{1-q} \, \varepsilon^{-2} \left(
    \sum_{i=0}^k \hat{c}_{3,0}^{(i)} c_{2,0} h_0^{\delta_i}
    + \sum_{i=0}^k \hat{c}_{3,0}^{(i)} \left(\frac{c_{2,0} c_2 c_3}{c_{3,0}}\right)^{1/2}
    L_{\alpha} h_0^{\delta_i} \right. \notag \\
    & + \hat{c}_3^{(0)} \left(\frac{c_2 c_{2,0} c_{3,0}}{c_3} \right)^{1/2}
    L_{\alpha} + \sum_{i=1}^k \hat{c}_3^{(i)} \left(\frac{c_2 c_{2,0} c_{3,0}}{c_3} \right)^{1/2}
    \frac{T^{\delta_i}-h_{L_{\alpha}}^{\delta_i}}{M^{\delta_i}-1} \notag \\
    & + \left. \hat{c}_3^{(0)} c_2 L_{\alpha}^2 + \sum_{i=1}^k \hat{c}_3^{(i)} c_2 L_{\alpha}
    \frac{T^{\delta_i}-h_{L_{\alpha}}^{\delta_i}}{M^{\delta_i}-1} \right) \notag \\
    \geq & \, \frac{T}{1-q} \, \varepsilon^{-2} \left(
    \sum_{i=0}^k \hat{c}_{3,0}^{(i)} c_{2,0} T^{\delta_i}
    + \left( \sum_{i=0}^k \hat{c}_{3,0}^{(i)} \left( \frac{c_{2,0} c_2 c_3}{c_{3,0}} \right)^{1/2}
    T^{\delta_i} \right. \right. \notag \\
    & + \left. \hat{c}_3^{(0)} \left(\frac{c_2 c_{2,0} c_{3,0}}{c_3} \right)^{1/2}
    + \sum_{i=1}^k \hat{c}_3^{(i)} c_2
    \frac{T^{\delta_i}-h_{L_{\alpha}}^{\delta_i}}{M^{\delta_i}-1} \right) \notag \\
    & \times \left( \frac{\log(\varepsilon^{-1})}{\alpha \log(M)} + \frac{\log(q^{-\frac{1}{2}}
    c_{1,\alpha} T^{\alpha})}{\alpha \log(M)} \right)
    + \hat{c}_3^{(0)} c_2 \left( \frac{\log(\varepsilon^{-1})}{\alpha \log(M)} \right)^2 \notag \\
    & + 2 \hat{c}_3^{(0)} c_2 \frac{\log(\varepsilon^{-1}) \log(q^{-\frac{1}{2}}
    c_{1,\alpha} T^{\alpha})}{\alpha^2 (\log(M))^2}
    + \hat{c}_3^{(0)} c_2 \frac{\big(\log(q^{-\frac{1}{2}} c_{1,\alpha} T^{\alpha})\big)^2}{(\alpha \log(M))^2} \notag \\
    & + \left. \sum_{i=1}^k \hat{c}_3^{(i)} \left(\frac{c_2 c_{2,0} c_{3,0}}{c_3} \right)^{1/2}
    \frac{T^{\delta_i}-h_{L_{\alpha}}^{\delta_i}}{M^{\delta_i}-1} \right)
    \label{Proof-Main-Prop-Lower-Bound-Eq-alg}
\end{align}
}
where $\hat{c}_{3,L_{\alpha}}^{(i)}=\hat{c}_{3}^{(i)}$,
$c_{2,L_{\alpha}}=c_2$, $c_{3,L_{\alpha}}=c_3$,
$\beta_{L_{\alpha}}=\beta$ and $\gamma_{L_{\alpha}}=\gamma$ for
$\hat{Y}_{ML(\alpha,\alpha)}$.

\emph{Upper bound for $\beta=\gamma$.}
Next, we calculate for $\beta=\gamma$ the upper bound
\begin{align}
    &C(\hat{Y}_{ML(\alpha,p)}) (\varepsilon)\notag\\
    \leq &\, \frac{T}{1-q} \, \varepsilon^{-2} \Bigg[
    \sum_{i=0}^k \hat{c}_{3,0}^{(i)} c_{2,0} T^{\delta_i}
    + \sum_{i=1}^k \hat{c}_3^{(i)} \left( \frac{c_{2,0} c_{3,0} c_2}{c_3} \right)^{1/2}
    \Lambda_i\notag \\
    & + \sum_{i=0}^k \hat{c}_{3,0}^{(i)} \left( \frac{c_{2,0} c_{2,L_p} c_{3,L_p}}{c_{3,0}} \right)^{1/2}
    T^{\delta_i} h_{L_p}^{\frac{\beta_{L_p}-\gamma_{L_p}}{2}} \notag \\
    & + \left( 
    \hat{c}_3^{(0)} \left( \frac{c_{2,0} c_{3,0} c_2}{c_3} \right)^{1/2}
    + \sum_{i=0}^k \hat{c}_{3,0}^{(i)} \left( \frac{c_{2,0} c_2 c_3}{c_{3,0}} \right)^{1/2}
    T^{\delta_i} \right. \notag \\
    & + \hat{c}_3^{(0)} \left( \frac{c_2 c_{2,L_p} c_{3,L_p}}{c_3} \right)^{1/2}
    h_{L_p}^{\frac{\beta_{L_p}-\gamma_{L_p}}{2}} \notag \\
    & + \left. \sum_{i=0}^k \hat{c}_{3,L_p}^{(i)} \left( \frac{c_2 c_3 c_{2,L_p}}{c_{3,L_p}} \right)^{1/2}
    h_{L_p}^{\frac{\beta_{L_p}-\gamma_{L_p}}{2}+\delta_i}
    + \sum_{i=1}^k \hat{c}_3^{(i)} c_2 \Lambda_i \right) \notag \\
    & \times \left( \frac{\log(\varepsilon^{-1})}{p \log(M)}
    + \frac{\log(q^{-\frac{1}{2}} c_{1,p} T^p)}{p \log(M)} \right) \notag \\
    & + \hat{c}_3^{(0)} c_2 \left( \frac{\log(\varepsilon^{-1})}{p \log(M)}
    + \frac{\log(q^{-\frac{1}{2}} c_{1,p} T^p)}{p \log(M)} \right)^2 \notag \\
    & + \sum_{i=1}^k \hat{c}_3^{(i)} \left( \frac{c_2 c_{2,L_p} c_{3,L_p}}{c_3} \right)^{1/2}
    \Lambda_i h_{L_p}^{\frac{\beta_{L_p}-\gamma_{L_p}}{2}} \notag \\
    & + \sum_{i=0}^k \hat{c}_{3,L_p}^{(i)} \left(
    \left( \frac{c_{2,0} c_{3,0} c_{2,L_p}}{c_{3,L_p}} \right)^{1/2}
    h_{L_p}^{\frac{\beta_{L_p}-\gamma_{L_p}}{2}+\delta_i} +
    c_{2,L_p} h_{L_p}^{\beta_{L_p}-\gamma_{L_p}+\delta_i} \right) \Bigg] \notag \\
    & + T \sum_{i=0}^k \left( \hat{c}_{3,0}^{(i)} T^{\delta_i-\gamma}
    + \hat{c}_3^{(i)} \frac{(M^{-1} T)^{\delta_i-\gamma} - h_{L_p}^{\delta_i-\gamma}}{1-M^{\gamma-\delta_i}}
    + \hat{c}_{3,L_p}^{(i)} h_{L_p}^{\delta_i-\gamma_{L_p}} \right)
    \label{Proof-Main-Prop-Upper-Bound-Eq-alg}
\end{align}
where we applied the relation (\ref{L-1-upper-bound}).

\emph{Proof of ii).}
Suppose that $\beta_{L_p} \geq \gamma_{L_p}$ and $\gamma,
\gamma_{L_p} \leq 2p$. Then, we get from the upper bound
(\ref{Proof-Main-Prop-Upper-Bound-Eq-alg}) that
\[C(\hat{Y}_{ML(\alpha,p)})(\varepsilon) = O( \varepsilon^{-2} (
\log(\varepsilon) )^2 ).\] Further, comparing the lower and the upper
bounds (\ref{Proof-Main-Prop-Lower-Bound-Eq-alg}) and
(\ref{Proof-Main-Prop-Upper-Bound-Eq-alg}), we asymptotically obtain
that
\begin{align}
    \lim_{\varepsilon \to 0} \frac{C(\hat{Y}_{ML(\alpha,\alpha)})(\varepsilon)}
    {C(\hat{Y}_{ML(\alpha,p)})(\varepsilon)}
    \geq \lim_{\varepsilon \to 0} \frac{\frac{T}{1-q} \varepsilon^{-2} \hat{c}_3^{(0)}
    c_2 \left( \frac{\log(\varepsilon^{-1})}{\alpha \, \log(M)} \right)^2
    + o(\varepsilon^{-2} (\log(\varepsilon))^2)}
    {\frac{T}{1-q} \varepsilon^{-2} \hat{c}_3^{(0)}
    c_2 \left( \frac{\log(\varepsilon^{-1})}{p \, \log(M)} \right)^2
    + o(\varepsilon^{-2} (\log(\varepsilon))^2)}
    = \frac{p^2}{\alpha^2}
\end{align}
which proves statement (\ref{Main-Prop-Improvement-Aussage2}).
This completes the proof.
\end{proof}
\begin{Bem}
    Especially, if $c_3 = \hat{c}_3
    $ and $c_{3,L_p}=\hat{c}_{3,L_p}
    $, then it follows in case of $\beta<\gamma$ and $\beta<2p$ that
    \begin{equation} \label{Main-Prop-Improvement-Aussage-Bemerkung-1}
        \lim_{\varepsilon \to 0} \frac{C(\hat{Y}_{ML(p,p)})(\varepsilon)}
        {C(\hat{Y}_{ML(\alpha,p)})(\varepsilon)} \geq M^{\gamma-\beta} \left(
        1 - M^{\frac{\beta-\gamma}{2}}
        \left(1 - \left( \frac{c_2 c_3}{c_{2,L_p} c_{3,L_p}} \right)^{1/2} \right)
        \right)^{-2} \, .
    \end{equation}
    Thus, if $c_2 c_3 < c_{2,L_p} c_{3,L_p}$ it follows directly that
    \begin{equation} \label{Main-Prop-Improvement-Aussage-Bemerkung-2}
        \lim_{\varepsilon \to 0} \frac{C(\hat{Y}_{ML(p,p)})(\varepsilon)}
        {C(\hat{Y}_{ML(\alpha,p)})(\varepsilon)} > 1 \, .
    \end{equation}
\end{Bem}
%
%
\section{Numerical examples in case of SDEs}
\label{Sec4:Numerical-Examples-SDEs}
%
For illustration of the improvement that can be realized with the
proposed mo\-di\-fied multi-level Monte Carlo estimator, we consider
the problem of weak approximation for stochastic differential
equations (SDEs)
\begin{equation} \label{SDE-general}
  \mathrm{d} X_t = a(X_t) \, \mathrm{d}t + \sum_{j=1}^m b^j(X_t) \, \mathrm{d}B_t^j
\end{equation}
with initial value $X_{t_0}=x_0 \in \mathbb{R}^d$ driven by
$m$-dimensional Brownian motion.

In the following, we compare for several numerical examples the root
mean-square errors \eqref{Sec2:root-mean-square-error} versus the
computational costs for the multi-level Monte Carlo estimator
$\hat{Y}_{ML}$ proposed in \cite{Gil08b,Gil09a,Gil08a} and described
in Section~\ref{Section2:MLMC-Simulation-Original} with the proposed
modified multi-level Monte Carlo estimator $\hat{Y}_{ML(\alpha,p)}$
described in Section~\ref{Sec3:Improved-MLMC-Estimator}.
As a measure for the computational costs, we count the number of
evaluations of the drift and diffusion functions taking into account
the dimension $d$ of the solution process as well as the dimension
$m$ of the driving Brownian motion.

In the following, we consider on each level $l=0,1, \ldots, L$ an
equidistant discretization $I_{h_l}=\{t_0, \ldots,
t_{\frac{T}{2^l}}\}$ of $[t_0,T]$ with step size
$h_l=\frac{T}{2^l}$. Further, we denote by $Y_n=Y_{t_n}$ the
approximation at time $t_n$. In case of the multi-level Monte Carlo
estimator $\hat{Y}_{ML}$ we apply on each level $l=0,1, \ldots, L$
the Euler-Maruyama scheme on the grid $I_{h_l}$ given by $Y_0=x_0$
and
\begin{equation}
  Y_{n+1} = Y_n + a(Y_n) \, h_n + \sum_{j=1}^m b^j(Y_n) \, I_{(j),n}
\end{equation}
where $h_n=h_l$ and $I_{(j),n} = B_{t_{n+1}}^j - B_{t_n}^j$ for
$n=0,1, \ldots, \tfrac{T}{2^l}-1$. The Euler-Maruyama scheme
converges with order $1/2$ in the mean-square sense and with order
$\alpha=1$ in the weak sense to the solution of the considered SDE
\eqref{SDE-general} at time $T$ \cite{KP99}.

On the other hand, for the modified multi-level Monte Carlo
estimator $\hat{Y}_{ML(\alpha,p)}$ the Euler-Maruyama scheme is
applied on levels $0,1, \ldots, L_p-1$ whereas on level $L_p$ a
second order weak stochastic Runge-Kutta (SRK) scheme RI6 proposed
in \cite{Roe09} is applied. The SRK scheme RI6 on level $L_p$ is
defined on the grid $I_{h_{L_p}}$ by $\check{Y}_0=x_0$,
\begin{equation}
  \begin{split}
  \check{Y}_{n+1} &= \check{Y}_n + \tfrac{1}{2} \left( a(\check{Y}_n) + a(\Upsilon) \right) \, h_{n}
  + \tfrac{1}{2} \sum_{k=1}^m \left( b^k(\Upsilon_+^{(k)}) - b^k(\Upsilon_-^{(k)}) \right)
  \, \tfrac{\hat{I}_{(k,k),n}}{\sqrt{h_{n}}} \\
  &+ \sum_{k=1}^m \left( \tfrac{1}{2} b^k(\check{Y}_n) + \tfrac{1}{4} b^k(\Upsilon_+^{(k)})
  + \tfrac{1}{4} b^k(\Upsilon_-^{(k)}) \right) \, I_{(k),n} \\
  &+ \tfrac{1}{2} \sum_{k=1}^m \left( b^k(\hat{\Upsilon}_+^{(k)}) - b^k(\hat{\Upsilon}_-^{(k)})
  \right) \, \sqrt{h_{n}} \\
  &- \sum_{k=1}^m \left(\tfrac{1}{2} b^k(\check{Y}_n) - \tfrac{1}{4} b^k(\hat{\Upsilon}_+^{(k)})
  - \tfrac{1}{4} b^k(\hat{\Upsilon}_-^{(k)}) \right) \, I_{(k),n}
  \end{split}
\end{equation}
where $h_n=h_{L_p}$ and $I_{(k),n} = B_{t_{n+1}}^k - B_{t_n}^k$ for
$n=0,1, \ldots, \tfrac{T}{2^{L_p}}-1$ with stages
\begin{equation}
  \begin{split}
  \Upsilon &= \check{Y}_n + a(\check{Y}_n) \, h_{n} + \sum_{j=1}^m b^j(\check{Y}_n) \, I_{(j),n} , \\
  \Upsilon_{\pm}^{(k)} &= \check{Y}_n + a(\check{Y}_n) \, h_{n} \pm b^k(\check{Y}_n) \, \sqrt{h_{n}} ,
  \quad \quad
  \hat{\Upsilon}_{\pm}^{(k)} = \check{Y}_n \pm
  \sum_{\substack{j=1 \\ j \neq k}}^m b^j(\check{Y}_n) \, \tfrac{\hat{I}_{(k,j),n}}{\sqrt{h_{n}}}
  \end{split}
\end{equation}
where $\hat{I}_{(k,k),n} = \tfrac{1}{2} (I_{(k),n}^2-h_{n})$ and
\begin{equation}
  \hat{I}_{(k,j)_n} = \begin{cases}
    \tfrac{1}{2} (I_{(k),n} I_{(j),n} - \sqrt{h_{n}} \tilde{I}_{(k),n}) & \text{ if } k<j \\
    \tfrac{1}{2} (I_{(k),n} I_{(j),n} + \sqrt{h_{n}} \tilde{I}_{(j),n}) & \text{ if } j<k
  \end{cases}
\end{equation}
based on independent random variables $\tilde{I}_{(k),n}$ with
$\Prob(\tilde{I}_{(k),n} = \pm \sqrt{h_{n}})=\tfrac{1}{2}$.
Thus, we have $\alpha=1$ and $p=2$ for the modified multi-level
Monte Carlo estimator $\hat{Y}_{ML(\alpha,p)}$ in the following.
Further, for both schemes the variance decays with the same order as
the computational costs increase, i.~e.\
$\beta=\beta_{L_p}=\gamma=\gamma_{L_p}=1$. Then, the optimal order
of convergence attained by the multi-level Monte Carlo method is
$O(\varepsilon^{-2} (\log(\varepsilon))^2)$ due to
Theorem~\ref{Main-Theorem-Giles}. For the presented simulations, we
denote by MLMC EM the numerical results for $\hat{Y}_{ML}$ based on
the Euler-Maruyama scheme only and by MLMC SRK the results for
$\hat{Y}_{ML(\alpha,p)}$ based on the combination of the
Euler-Maruyama scheme and the SRK scheme RI6.

\begin{figure}
\begin{center}
\includegraphics[width=5.5cm]{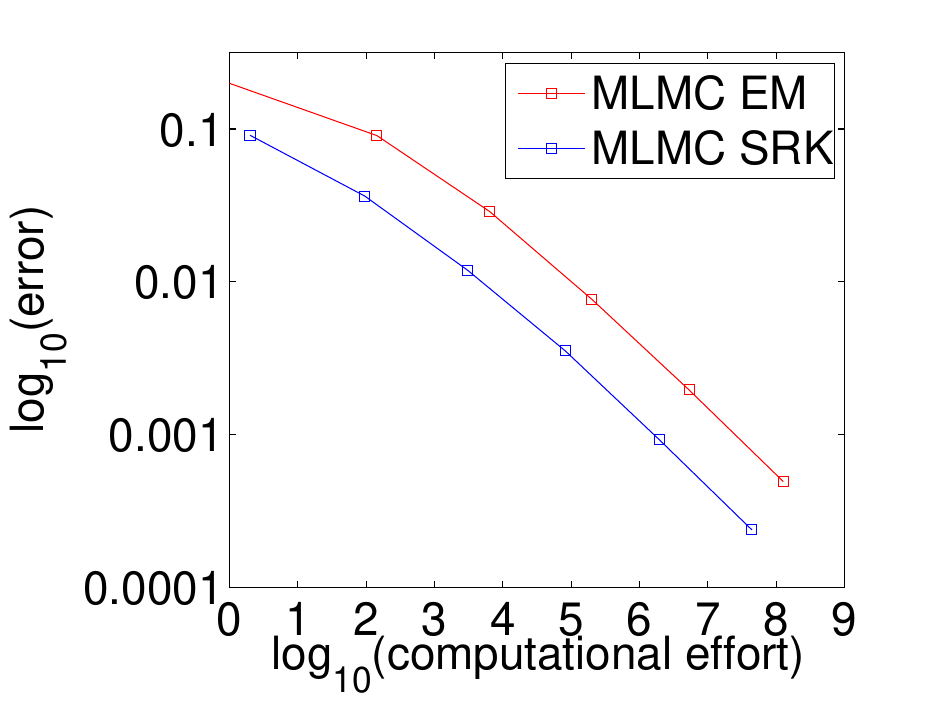}
\qquad
\includegraphics[width=5.5cm]{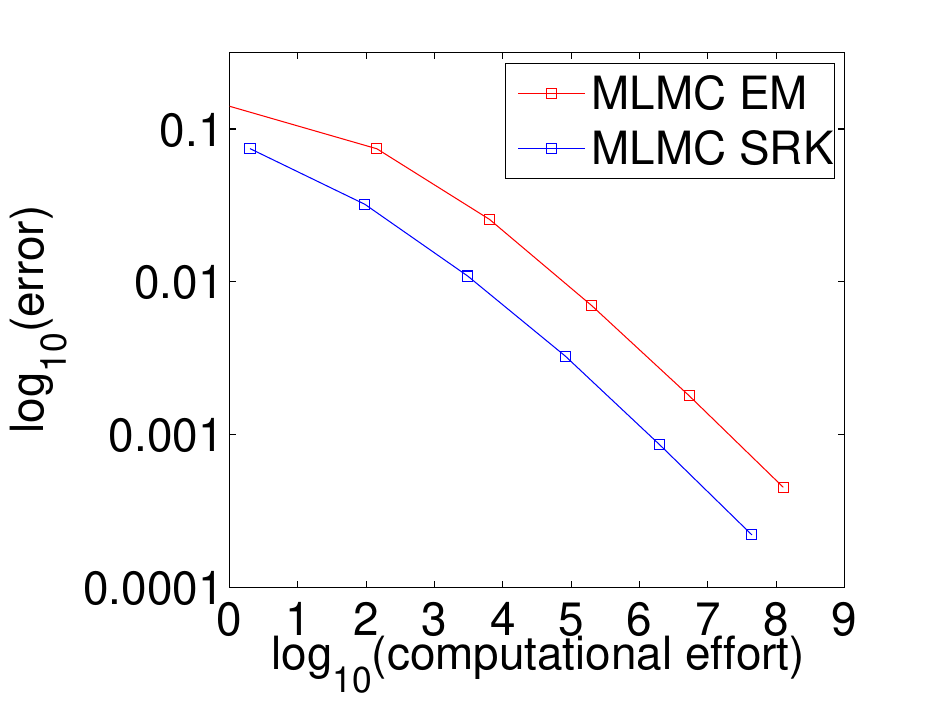}
\end{center}
\caption{Error vs.\ computational effort for SDE~\eqref{Test-SDE-1} using $f(x)=x$ (left) and
$f(x)=x^2$ (right).}
\label{Bild-test1}
\end{figure}

As a first example, we consider the scalar linear SDE with $d=m=1$
given by
\begin{equation} \label{Test-SDE-1}
        {\mathrm{d}} X_t = r \, X_t \, {\mathrm{d}}t
        + \sigma X_t \, {\mathrm{d}}B_t \, ,
        \quad X_0 = 0.1 \, ,
\end{equation}
using the parameters $r=1.5$ and $\sigma=0.1$. We choose $T=1$ and
apply the functionals $f(x)=x$ and $f(x)=x^2$, see
Figure~\ref{Bild-test1}. The presented simulations are calculated
using the prescribed error bounds $\varepsilon=4^{-j}$ for $j=0,1,
\ldots, 5$. In Figure~\ref{Bild-test1} we can see the significantly
reduced computational effort for the estimator $\hat{Y}_{ML(1,2)}$
(MLMC SRK) compared to the estimator $\hat{Y}_{ML}$ (MLMC EM) in
case of a linear and a nonlinear functional.

\begin{figure}
\begin{center}
    \includegraphics[width=5.5cm]{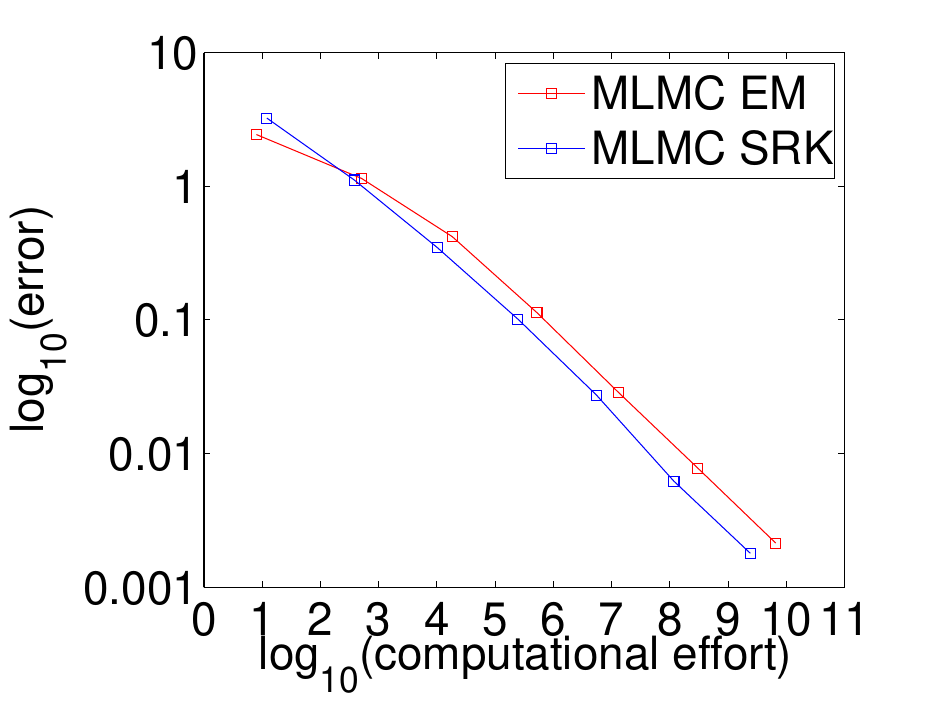}
    \qquad
    \includegraphics[width=5.5cm]{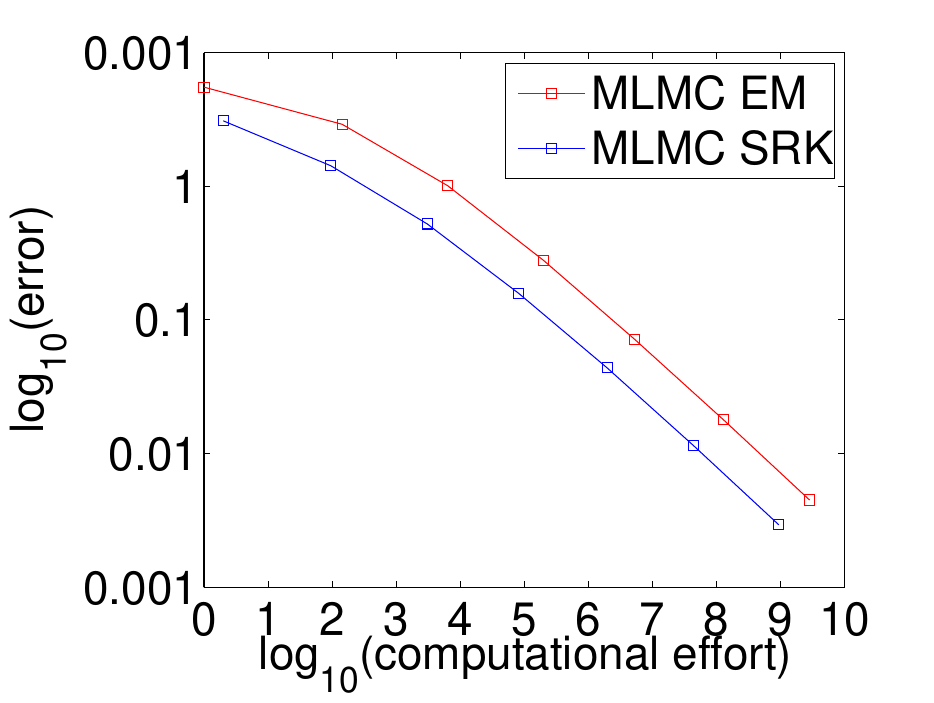}
\end{center}
\caption{Error vs.\ computational effort for the nonlinear SDE~\eqref{Test-SDE-2} (left) and
SDE~\eqref{Test-SDE-3} (right) with non-commutative noise.}
\label{Bild-test2}
\end{figure}

The second example is a nonlinear scalar SDE with $d=m=1$ given by
\begin{equation} \label{Test-SDE-2}
        {\mathrm{d}} X_t = \tfrac{1}{2} X_t +\sqrt{X_t^2+1} \, {\mathrm{d}}t
        + \sqrt{X_t^2+1} \, {\mathrm{d}}B_t \, ,
        \quad X_0 = 0 \, .
\end{equation}
We apply the functional
\[f(x)= (\log(x + \sqrt{x^2 + 1}))^3 - 6
(\log(x + \sqrt{x^2 + 1}))^2 + 8 \log(x + \sqrt{x^2 + 1}).\] Then,
the approximated expectation is given by \[\Erw(f(X_t)) = t^3-3 t^2 + 2
t.\] Here, the results presented in Figure~\ref{Bild-test2} (left)
are calculated for $T=2$ applying the prescribed error bounds
$\varepsilon=4^{-j}$ for $j=0,1, \ldots, 6$. Here, the improved
estimator $\hat{Y}_{ML(1,2)}$ performs much better than
$\hat{Y}_{ML}$ also for nonlinear functionals and a nonlinear SDE.
Finally, we consider a nonlinear multi-dimensional SDE with a $d=4$
dimensional solution process driven by an $m=6$ dimensional Brownian
motion with non-commutative noise:
\begin{multline} \label{Test-SDE-3}
    {\mathrm{d}} \begin{pmatrix} X_t^1 \\ X_t^2 \\ X_t^3 \\ X_t^4 \end{pmatrix}
     = \begin{pmatrix} \frac{243}{154} X_t^1 - \frac{27}{77} X_t^2 + \frac{23}{154} X_t^3
     - \frac{65}{154} X_t^4 \\
    \frac{27}{77} X_t^1 - \frac{243}{154} X_t^2 + \frac{65}{154}
    X_t^3 - \frac{23}{154} X_t^4 \\
    \frac{5}{154} X_t^1 - \frac{61}{154} X_t^2 + \frac{162}{77}
    X_t^3 - \frac{36}{77} X_t^4 \\
    \frac{61}{154} X_t^1 - \frac{5}{154} X_t^2 + \frac{36}{77}
    X_t^3 - \frac{162}{77} X_t^4 \end{pmatrix} \, {\mathrm{d}}t \\
    + \frac{1}{9} \sqrt{(X_t^2)^2 + (X_t^3)^2 + \frac{2}{23}}
    \begin{pmatrix} \tfrac{1}{13} \\ \tfrac{1}{14} \\ \tfrac{1}{13} \\
    \tfrac{1}{15} \end{pmatrix} \, {\mathrm{d}} B_t^1
    + \frac{1}{8} \sqrt{(X_t^4)^2 + (X_t^1)^2 + \frac{1}{11}}
    \begin{pmatrix} \frac{1}{14} \\ \frac{1}{16} \\ \frac{1}{16} \\
    \frac{1}{12} \end{pmatrix} \, {\mathrm{d}} B_t^2 \\
    + \frac{1}{12} \sqrt{(X_t^1)^2 + (X_t^2)^2 + \frac{1}{9}}
    \begin{pmatrix} \frac{1}{6} \\ \frac{1}{5} \\ \frac{1}{5} \\
    \frac{1}{6} \end{pmatrix} \, {\mathrm{d}} B_t^3
    + \frac{1}{14} \sqrt{(X_t^3)^2 + (X_t^4)^2 + \frac{3}{29}}
    \begin{pmatrix} \frac{1}{8} \\ \frac{1}{9} \\ \frac{1}{8} \\
    \frac{1}{9} \end{pmatrix} \, {\mathrm{d}} B_t^4 \\
    + \frac{1}{10} \sqrt{(X_t^1)^2 + (X_t^3)^2 + \frac{1}{13}}
    \begin{pmatrix} \frac{1}{11} \\ \frac{1}{15} \\ \frac{1}{13} \\
    \frac{1}{11} \end{pmatrix} \, {\mathrm{d}} B_t^5
    + \frac{1}{11} \sqrt{(X_t^2)^2 + (X_t^4)^2 + \frac{2}{25}}
    \begin{pmatrix} \frac{1}{12} \\ \frac{1}{13} \\ \frac{1}{16} \\
    \frac{1}{13} \end{pmatrix} \, {\mathrm{d}} B_t^6
\end{multline}
\noindent
with initial condition $X_0 = (\tfrac{1}{8}, \tfrac{1}{8}, 1,
\tfrac{1}{8})^T$. Then, the approximated first moment of the
solution is given by $E(X_T^i)= X_0^i \, \exp(2T)$ for $i=1,2,3,4$.
The simulation results calculated at $T=1$ for the error bounds
$\varepsilon=4^{-j}$ for $j=0,1, \ldots, 6$ are presented in
Figure~\ref{Bild-test2} (right). Again, in the multi-dimensional
non-commutative noise case the proposed estimator
$\hat{Y}_{ML(1,2)}$ needs significantly less computational effort
compared to the estimator $\hat{Y}_{ML}$ which reveals the
theoretical results \eqref{Main-Prop-Improvement-Aussage2} in
Proposition~\ref{Main-Prop-Improvement}.
%
%
%
\section{Conclusions}
\label{Sec5:Conclusions}
%
%
In this paper we proposed a modification of the multi-level Monte
Carlo method introduced by M.\ Giles which combines approximation
methods of different orders of weak convergence. This modified
multi-level Monte Carlo method attains the same mean square order of
convergence like the originally proposed method that is in some
sense optimal. However, the newly proposed multi-level Monte Carlo
estimator can attain significantly reduced computational costs. As
an example, there is a reduction of costs by a factor $(p/\alpha)^2$
for the problem of weak approximation for SDEs driven by Brownian
motion in case of $\beta=\gamma$. This has been approved by some
numerical examples for the case of $p=2$ and $\alpha=1$ where four
times less calculations are needed compared to the standard
multi-level Monte Carlo estimator. Here, we want to point out that
there also exist higher order weak approximation schemes, e.~g.\
$p=3$ in case of SDEs with additive noise \cite{De10}, that may
further improve the benefit of the modified multi-level Monte Carlo
estimator. Future research will consider the application of this
approach to, e.g., more general SDEs like SDEs driven by L\'{e}vy
processes \cite{Der11} or fractional Brownian motion \cite{KNP11}
and to the numerical solution of SPDEs \cite{SchGit11}. Further, the
focus will be on numerical schemes that feature not only high orders
of convergence but also minimized constants for the variance
estimates.

\end{document}